\documentclass[12pt]{amsart}%
\usepackage{graphicx}
\usepackage{amssymb}
\usepackage{amsmath}
\usepackage{amsfonts}
\numberwithin{equation}{section}
\begin{document}
\title[Delayed heat equation]{A note on exponential-type solutions for the linear, delayed heat partial
differential equation }
\author[I.H. Herron and R. E. Mickens]{Isom H. Herron \\
Department of Mathematical Sciences \\
Rensselaer Polytechnic Institute \\
Troy, NY 12180 USA \\
Email: herroi@rpi.edu \\
 and \\
Ronald E. Mickens \\
Department of Physics\\
Clark Atlanta University\\
Atlanta, GA 30314 USA \\
Email: rmickens@cau.edu}
\begin{abstract}
We construct a class of exponential type solutions for the linear, delayed
heat equation. These representations may be used to provide a priori ansatzes
for certain boundary and/or initial-value problems arising in heat transfer.
Several of the important mathematical properties of the representations are
examined, including a discussion of the dependence on the delay parameter.

\end{abstract}
\maketitle


{
{

\section{Introduction}

The purpose of this Note is to construct a special class of exact solutions to
the linear, delayed heat partial differential equation (PDE)
(\cite{Jordan2007}, \cite{Jordan/Dai/Mickens}, \cite{Joseph&Preziosi},
\cite{Tzou})%

\begin{equation}
u_{t}\left(  x,t+\tau\right)  =Du_{xx}\left(  x,t\right)  , \label{(1.1)}%
\end{equation}
where $\left(  x,t\right)  $ are respectively space and time variables,
$\ \tau\geq0$ is time delay, $D>0$. The importance of this equation is
indicated not only by its application to a broad range of phenomena in the
physical and engineering sciences, but also the number of papers devoted to
investigating the nature of its solutions; see the above indicated references,
along with that of Tzou \cite{Tzou} and the works cited within them.

Our major result is that the delayed heat equation (DHE) has a special class
of solutions which take the form
\begin{equation}
u\left(  x,t,\tau,a,b\right)  =e^{ax+bt}, \label{(1.2)}%
\end{equation}
where $a$ depends on $b,\ $i.e.%
\begin{equation}
a=a\left(  \tau,D,b\right)  . \label{(1.3)}%
\end{equation}
Further since the DHE is linear, general solutions may be constructed by
linear combination, i.e.
\begin{equation}
u\left(  x,t,\tau\right)  ={\sum\kern-1.1em\int}f\left(  b\right)  e^{a\left(
b+bt\right)  }db \label{1.4}%
\end{equation}
where $f\left(  b\right)  $ is such that this mathematical expression is
defined. Note that the explicit dependence on $D$ has been dropped and the
symbol in front on the right side of (\ref{1.4}) indicates \textquotedblleft
integration\textquotedblright\ over $b$ for continuous values of $b\ $and
\textquotedblleft sums\textquotedblright\ for discrete values of $b.$
Generally, which occurs depends on the particular initial- and/or
boundary-values considered for a given problem.

Of interest is that for boundary-value problems, the so-called Lambert
$W\ $function (\cite{Corless/etal}, \cite{Valluri}) makes an appearance.
Previously an instability was uncovered in the analysis of the heat equation
with time delay when solved on a finite spatial interval
\cite{Jordan/Dai/Mickens}, by making use of the $W-$function.

The next section is devoted to the explicit construction of these exponential
solutions. This is followed by an examination of the mathematical structure of
the solutions, along with a brief discussion of some possible next steps.

\section{Exponential solutions}

From (\ref{(1.2)}) it follows that%
\begin{align}
u\left(  x,t+\tau\right)   &  =e^{b\tau}\left(  ...\right)  ,\label{2.1a}\\
u_{t}\left(  x,t+\tau\right)   &  =be^{b\tau}\left(  ...\right)  ,
\label{2.1b}\\
u_{xx}(x,t)  &  =a^{2}\left(  ...\right)  \label{2.1c}%
\end{align}
where%
\begin{equation}
\left(  ...\right)  =e^{ax+bt}. \label{2.1d}%
\end{equation}
Therefore, substituting these expressions into eqn. (\ref{(1.1)}) and
cancelling a common factor gives%
\begin{equation}
be^{b\tau}=Da^{2}, \label{2.2}%
\end{equation}%
\begin{equation}
a^{2}=\left(  \frac{1}{D}\right)  be^{b\tau}. \label{2.3}%
\end{equation}
There are two nontrivial cases to consider.\newline$b>0:$\newline This case
gives%
\begin{equation}
a=\left(  \pm\right)  a_{+}=\left(  \pm\right)  \left(  \sqrt{\frac{b}{D}%
}\right)  e^{b\tau/2} \label{2.4}%
\end{equation}
$b<0:$\newline For this situation, we have%
\begin{equation}
a=\left(  \pm\right)  ia_{-}=\left(  \pm\right)  i\left(  \sqrt{\frac
{\left\vert b\right\vert }{D}}\right)  e^{-\left\vert b\right\vert \tau/2},
\label{2.5}%
\end{equation}
where $i=\sqrt{-1}$ .

Making use of the Euler relations%
\begin{align}
\cos\theta &  =\frac{e^{i\theta}+e^{-i\theta}}{2},\sin\theta=\frac{e^{i\theta
}-e^{-i\theta}}{2i},\label{2.6a}\\
\cosh\theta &  =\frac{e^{\theta}+e^{-\theta}}{2},\sinh\theta=\frac{e^{\theta
}-e^{-\theta}}{2}, \label{2.6b}%
\end{align}
we obtain, after some algebraic manipulations, the results%
\begin{align}
u\left(  x,t,\tau,b>0\right)   &  \equiv u_{+}\left(  x,t,\tau,b\right)
,\label{2.7a}\\
u\left(  x,t,\tau,b<0\right)   &  \equiv u_{-}\left(  x,t,\tau,b\right)  .
\label{2.7b}%
\end{align}
In detailed, explicit form, we have%
\begin{equation}
u_{+}\left(  x,t,\tau,b\right)  =e^{bt}\left[  A_{1}\left(  b\right)
\cosh\left(  a_{+}x\right)  +A_{2}\left(  b\right)  \sinh\left(
a_{+}x\right)  \right]  \label{2.8}%
\end{equation}
and%
\begin{equation}
u_{-}\left(  x,t,\tau,b\right)  =e^{-\left\vert b\right\vert t}\left[
B_{1}\left(  b\right)  \cos\left(  a_{-}x\right)  +B_{2}\left(  b\right)
\sin\left(  a_{-}x\right)  \right]  , \label{2.9}%
\end{equation}
where $\left(  a_{+},a_{-}\right)  $ are given respectively, by (\ref{2.4})
and (\ref{2.5}) and $A_{1},A_{2},B_{1},B_{2}$ are arbitrary functions of $b.$

Using the fact that the DHE is linear \ allows for the construction of more
general solutions through the formulation of linear combination of solutions.
Doing this gives the result%
\begin{equation}
u\left(  x,t,\tau\right)  ={\sum_{b>0}\kern-1.1em\int}u_{+}\left(
x,t,\tau,b\right)  db+{\sum_{b<0}\kern-1.1em\int}u_{-}\left(  x,t,\tau
,b\right)  db. \label{2.10}%
\end{equation}

Inspection of eqn.\thinspace(\ref{2.10}) indicates that the first term on the
right-side, gives rise to solutions unbounded in time, if given finite initial
and boundary conditions at $t=0.$ Thus, for bounded solutions, only the second
term on the right side should be considered.

\section{Discussion}

First, observe that from eqns.\thinspace(\ref{2.8}), (\ref{2.9}) and
(\ref{2.10}) that
\begin{equation}
\lim_{\tau\rightarrow0}u(x,t,\tau)=w\left(  x,t\right)  , \label{3.1}%
\end{equation}
where $w\left(  x,t\right)  $ is a solution to%
\begin{equation}
w_{t}\left(  x,t\right)  =Dw_{xx}\left(  x,t\right)  , \label{3.2}%
\end{equation}
which is the standard heat equation without delay. In this limit,
$\tau\rightarrow0$%
\begin{align}
a_{+}\left(  b,\tau\right)   &  \rightarrow\sqrt{\frac{b}{D}},b>0;
\label{3.3a}\\
a_{-}\left(  b,\tau\right)   &  \rightarrow\sqrt{\frac{\left\vert b\right\vert
}{D}},b<0. \label{3.3b}%
\end{align}
Consequently, the exponential type solutions of eqn. (\ref{(1.1)}) reduce to
the correct, corresponding solutions of (\ref{3.2}).

Second, for the initial-value/boundary-value problem $\left(  0\leq x\leq
L\right)  $%
\begin{align}
u\left(  x,t,\tau\right)   &  =f\left(  x\right)  =\text{given, }-\tau\leq
t\leq0\label{3.4a}\\
u\left(  0,t,\tau\right)   &  =0,\ \ \ u\left(  L,t,\tau\right)  =0,\ \ t>0,
\label{3.4b}%
\end{align}
the following condition is obtained%
\begin{equation}
a_{-}\left(  b,\tau\right)  =\frac{n\pi}{L},(n=1,2,3,...), \label{3.5}%
\end{equation}
and this can be rewritten to the form%
\begin{equation}
\left(  \frac{\left\vert b\right\vert }{D}\right)  e^{-\left\vert b\right\vert
\tau}=\frac{n^{2}\pi^{2}}{L^{2}} \label{3.6}%
\end{equation}

Let $z=-\left\vert b\right\vert \tau$, then we have%
\begin{equation}
ze^{z}=-\tau D\left(  \frac{n\pi}{L}\right)  ^{2}. \label{3.7}%
\end{equation}
Note that this equation may be written in the form of the Lambert $W$ function
\cite{Corless/etal}.\newline However, without the use of this function, we
have from a traditional knowledge of transcendental functions that all the
roots of
\[
ze^{z}+q=0,\text{\ for\ }q\text{ }\operatorname{real},
\]
have negative real parts if and only if
\[
0\,<q<\frac{\pi}{2}.
\]
Applied to this case this means%
\[
\tau D\left(  \frac{n\pi}{L}\right)  ^{2}<\frac{\pi}{2}.
\]
This inequality will not hold for $n$ sufficiently large.

\paragraph{Comment-1}

The extensive work of Pedro Jordan clearly indicates the perils of using
eqn.\thinspace(\ref{(1.1)}) as a mathematical model for heat transfer; see
\cite{Jordan/Dai/Mickens}. Further, the exponential solution is related to his
work. For example the solution given by (\ref{(1.2)})%

\[
u\left(  x,t,\tau,a,b\right)  =e^{ax+bt},
\]
is an alternative to Jordan's separation of variables ansatz
\cite{Jordan/Dai/Mickens}
\begin{equation}
u\left(  x,t,\tau,a,b\right)  =X\left(  x\right)  T\left(  t\right)  .
\label{3.8}%
\end{equation}
The major difference of the two procedures is in the appearance of the Lambert
$W-$function, now occurring in calculating solutions for $T\left(  t\right)
.$ This illustrates the fact that exponential and separation of variables
methods are related to each other. In particular the exponential ansatz is a
special case of the separation of variables method.

\paragraph{Comment-2}

Using a somewhat different methodology, Ismagilov et al \cite{Ismagiilov} came
to same conclusion as Jordan et al \cite{Jordan/Dai/Mickens} regarding the
stability of solutions for the heat equation with delay. In \cite{Ismagiilov},
a similar conclusion was reached for the wave equation with delay%
\begin{equation}
u_{tt}\left(  x,t+\tau\right)  =u_{xx}\left(  x,t\right)  . \label{3.9}%
\end{equation}
See also the paper \cite{Rodrigues}.

\paragraph{Comment-3}

Polynain and Zhurov \cite{Polyanin} have constructed a set of procedures for
calculating solutions to non-linear, delay, reaction-diffusion partial
differential equations. It might prove of value to see if this methodology can
also be usefully applied to the linear, delayed heat equation. Similarly, one
could apply the generating function technique (GFT) created by Robert Jackson
\cite{Jackson} to this equation to see if new types of solutions can be constructed.

\paragraph{Comment-4}

The above methodology can also be applied to the linear delayed
advection-diffusion PDE%
\begin{equation}
u_{t}\left(  x,t+\tau\right)  +\varepsilon u_{x}\left(  x,t\right)
=Du_{xx}\left(  x,t\right)  , \label{3.10}%
\end{equation}
where$\ \left(  \varepsilon,D\right)  $ are non-negative parameters.\newline
For this case, we have%
\begin{equation}
u\left(  x,t\right)  =e^{ax+bt} \label{3.11}%
\end{equation}
with%
\begin{equation}
Da^{2}-\varepsilon a-be^{b\tau}=0. \label{3.12}%
\end{equation}
Solving for $a$ gives%
\begin{equation}
a_{\pm}=\frac{1}{2D}\left[  \varepsilon\pm\sqrt{\varepsilon^{2}+4bDe^{b\tau}%
}\right]  \label{3.13}%
\end{equation}
where
\begin{equation}
a_{-}<0<a_{+},\ a_{+}>\left\vert a_{-}\right\vert . \label{3.14}%
\end{equation}
Therefore,%
\begin{equation}
u\left(  x,t\right)  ={\sum\kern-1.1em\int}\hat{u}\left(  b,x,t\right)  db
\label{3.15}%
\end{equation}
with%
\begin{equation}
\hat{u}\left(  b,x,t\right)  =e^{bt}\left[  A_{1}\left(  b\right)  e^{a_{+}%
x}+A_{2}\left(  b\right)  e^{-\left\vert a_{-}\right\vert x}\right]  ,
\label{3.16}%
\end{equation}
where $A_{1}\left(  b\right)  $ and $A_{2}\left(  b\right)  $ are functions of
$b $ such that for particular application, eqn. (3.15) exists.

Finally it should be mentioned that the major reason for the interest in the
delayed heat equation (DHE) is that it provides a model for heat transfer for
a broad range of phenomena. However, there are a variety of investigations
indicating that this model and approximation to it may not be valid
(\cite{Jordan2007}, \cite{Jordan/Dai/Mickens}, \cite{Joseph&Preziosi},
\cite{Tzou}). A major problem is that the microscopic, i.e. atomic structure
of matter, is not included in the mathematical formulation. How to overcome
this difficulty is a goal to fulfill by extending the results of the current
research efforts.

\textbf{Acknowledgement }REM thanks Dr.\thinspace Pedro M.\thinspace Jordan
(Naval Research Laboratory, Stennis Space Center, MS) for many fruitful and
deep discussions on the subject of heat transf er and the fundamental need for
an enhanced heat equation going beyond the standard delayed flux relation.

\end{document}